\documentclass[10pt]{ip-journal}
\usepackage{amssymb}
\usepackage{amsfonts}
\usepackage{graphicx}
\usepackage{amsmath}
\usepackage{color}

\newtheorem{theorem}{Theorem}[section]
\renewcommand{\theequation}{\thesection.\arabic{equation}}

\newtheorem{assumption}[theorem]{Assumption}

\newtheorem{lemma}[theorem]{Lemma}

\newtheorem{proposition}[theorem]{Proposition}
\newtheorem{remark}[theorem]{Remark}

\newcommand{\R}{\mathbb{R}}

\begin{document}

\title{PERTURBATION OF A GLOBALLY STABLE EQUILIBRIUM: APPLICATION ON AN AGE-STRUCTURED MODEL
$^{*}$\footnote{$^{*}$ This work was supported by the government of Canada's International Development Research Centre (IDRC), and within the framework of the AIMS Research for Africa Project No SNMCM2013014S. The author acknowledges also a Post AIMS-Senegal Grant No 001CORE140450006. This work is a chapter of the PhD work of YTK (corresponding author).}} 

\author{Yannick Tchaptchie Kouakep$^{**}$\footnote{$^{**}$ University of Ngaoundere, ERMIA, PO Box 454 Ndang, Ngaoundere, Cameroon. 
(e-mail: kouakep@aims-senegal.org). With LYCLAMO (PO Box 46, Ngaoundere, Cameroon) And AIMS - Cameroon, PO Box 608, Limbe, Cameroon.}}

\maketitle

\begin{abstract}
In this manuscript we consider an age structured epidemic system modelling the dynamics of transmission of immunizing disease like Hepatitis B virus. Our model takes into account age as well as two classes of infected individuals (chronic carriers and acute infected human). Based on the low infectivity of chronic carriers, we study the asymptotic behaviour of the system and, under some suitable assumptions, we prove the global stability of the endemic equilibrium point using perturbation arguments of semiflow. The intuitive conservation of the global stability under perturbation, needs in fact very technical tools and complex results.
\end{abstract}

\noindent {\bf AMS subject classification: 35Q92, 34K20, 92D30.}\\
\noindent {\bf Keywords: Age-structured model, perturbation arguments, asymptotic behaviour, Hepatitis B}

\pagestyle{myheadings}
\thispagestyle{plain}
\markboth{Yannick TCHAPTCHIE KOUAKEP}{PERTURBATION OF A GLOBALLY STABLE EQUILIBRIUM}

\section{Introduction}
We consider two age-structured systems modelling an immunizing disease like hepatitis B and link them to the main partially ''aggregated'' model to be studied in this manuscript. The first one uses the variables: 
\begin{itemize}
\item[a-] $s(t,a)$ the density of susceptible(s) at time $t$ with chronological age $a$;
\item[b-] $i(t,\tau)$ the density of infective(s) that will develop acute disease at time $t$ contaminated since a time $\tau$ ;
\item[c-] $j(t,\tau)$ the density of infective(s) that will not develop acute disease (asymptomatic carrier(s)) at time $t$ contaminated since a time $\tau$ (sometimes referred to as ``age of infection'' which plays a significant role in infective compartments $i$ and $j$)
\end{itemize}
and reads as
\begin{equation}\label{eq1.1}
\begin{split}
&(\partial_t+\partial_a)s(t,a)=-\mu s(t,a)-\lambda_0(t)s(t,a),\;\;t>0,\;a>0,\\
&s(t,0)=\Lambda,\\
&(\partial_t+\partial_\tau)i(t,\tau)=-(\mu +\gamma_i)i(t,\tau),\;t>0,\;\tau>0,\\
&(\partial_t+\partial_\tau)j(t,\tau)=-\left(\mu+\gamma_e\right) j(t,\tau),\;t>0,\;\tau>0,\\
&i(t,0)=\lambda_0(t)\int_0^\infty p_I(a)s(t,a)da,\;\;j(t,0)=\lambda_0(t)\int_0^\infty p_J(a)s(t,a)da.\\
\end{split}
\end{equation}

Here $\Lambda>0$ is some constant entering influx, $\mu>0$ is the natural death rate, $\gamma_i$ and $\gamma_e$ are the additional death rates due to the disease. We set also (for \eqref{eq1.1}-\eqref{eq1.2i}) the positive constants: $\nu_I=(\mu +\gamma_i)$ and $\nu_J=(\mu +\gamma_j)$. In addition $p_I\in L_+^\infty(0,\infty)$ is a given function such that $0\leq p_I(a)\leq 1$ $a.e.$ while $p_J(a)\equiv 1-p_I(a)$. Function $p_J$ represents the age-specific probability to become a chronic carrier when becoming infected at age $a$. Function $p_I$ denotes the probability to develop an acute infection when getting the infection at age $a$ (\cite{hw,kiire}). We refer to Edmunds et al.\cite{Ed93} for more explanation on the age-dependence susceptibility to the infection.

The force of infection $\lambda_0(t)$ (for \eqref{eq1.1}-\eqref{eq1.2i}) is expressed as
\begin{equation}\label{eq1.2}
\lambda_0(t)=\int_0^\infty \left( \beta_i(\tau) i(t,\tau)+\beta_j(\tau) j(t,\tau)\right)d\tau.
\end{equation}
Here in \eqref{eq1.2}, $\beta_i(\tau)$ and $\beta_j(\tau)$ denote the contact transmission rates between acute infected or chronic carriers with age of infection $\tau$ with susceptible individuals respectively.

Finally this first model is supplemented together with some initial data
\begin{equation}\label{eq1.2i}
\begin{split}
&s(0,.)=s_0(.)\in L^1_+(0,\infty)\\
&i(0,.)=i_0(.),\;\;j(0,.)=j_0(.)\text{ with }(i_0,j_0)\in \left[L^1_+(0,\infty)\right]^2.
\end{split}
\end{equation}

The second model with chronological age is
\begin{equation}\label{eq0}
\begin{cases}
\left[\partial_t+\partial_a+\mu\right]s(t,a)=-\lambda_1(t,a)s(t,a),\\
\left[\partial_t+\partial_a +\left(\mu_I+\mu+\omega_I\right)\right]i(t,a)=\lambda_1(t,a)p_I(a)s(t,a),\\
\left[\partial_t+\partial_a+\nu_J\right]j(t,a)=\lambda_1(t,a)p_J(a)s(t,a),\\
\left[\partial_t+\partial_a +\mu\right]r(t,a)=\mu_I i(t,a),
\end{cases}
\end{equation}
posed for $t>0$ and $a>0$, wherein $t$ denotes time while $a$ denotes the chronological age of individuals, that is the time since birth. Here $s(t,a)$ denotes the age-specific density of susceptible, $j(t,a)$ and $i(t,a)$ denotes respectively the the age-specific density of chronic carriers and acute infected individuals (that can be symptomatic or asymptomatic) for Hepatitis B virus (HBV) while $r(t,a)$ denotes the density of recovered individuals from HBV acute infection.

Parameter $\mu>0$ denotes the natural death rate, $\nu_I:=(\mu_I+\mu+\omega_I)>0$ and $\nu_J$ denote (for \eqref{eq0}-\eqref{eq000}) the exit rates associated to each infected class $i$ and $j$. The term $\nu_I$ gathers immunisation rate $\mu_I\geq 0$ due to acute infection, natural death rate $\mu$ and a possible additional death $\omega_I$ due to the infection; while $\nu_J$ corresponds to death of chronic carriers, that is natural death rate and an additional death rate due to chronic infection and its consequences. Here we leave out possible recovery from chronic disease. The term $\lambda_1(t,a)$ corresponds to the age-specific force of infection and follows the usual law of mass-action, that reads as
\begin{equation}\label{lam1}
\lambda_1(t,a)=\int_0^\infty \left[\beta_i(a,a')i(t,a')+\beta_j(a,a')j(t,a')\right]da'.
\end{equation}
Here in \eqref{lam1}, $\beta_i(a,a')$ and $\beta_j(a,a')$ denote the contact transmission rates between acute infected or chronic carriers of chronological age $a'$ with susceptible of chronological age $a$ respectively.

This problem \ref{eq0} is supplemented  together with the boundary conditions:
\begin{equation}\label{eq00}
\begin{cases}
s(t,0)=\Lambda,\text{ (constant influx)}\\
i(t,0)=j(t,0)=0,\text{ (no vertical transmission)},\\
r(t,0)=0,\text{ (no immunity at birth)},
\end{cases}
\end{equation}
and initial data
\begin{equation}\label{eq000}
s(0,a)=s_0(a),\;\;i(0,a)=i_0(a),\;\;j(0,a)=j_0(a),\;\;r(0,a)=r_0(a).
\end{equation}
Note that the $r-$components of the system decouples from the other and have therefore no impact upon the long time behaviour of the system. It will be omitted in the sequel.
Here recall that the boundary conditions $i(t,0)=j(t,0)=0$ correspond to omitting vertical transmission (see \cite{abbas,majid,m1,hw,kiire}).

In the absence of the disease, the evolution of the population density follows the following simple (chronological) age-structured equation:
\begin{equation*}
\left(\partial_t+\partial_a+\mu\right)s(t,a)=0,\;\;s(t,0)=\Lambda,
\end{equation*}
where $\mu>0$ denotes the natural death rate while $\Lambda>0$ corresponds to the birth rate, that is assumed to be a constant influx.

With 350 million people chronically infected, the Hepatitis B virus (HBV) pandemic is a global public health problem and constitutes one of the main cause of cirrhosis and hepatocellular carcinoma with a high rate of morbidity and mortality especially in South-East Asia and Sub-Saharan Africa struggling with a prevalence of HbsAg (Hepatitis B surface antigens) chronic carriers greater than 8$\%$ \cite{f1,1et,2et,Tiollais}. 
The class of HBV chronic carriers has a high risk of developing liver disease such as the one mentioned above.  

In highly endemic areas the epidemiology of HBV exhibits two main features: a low average age at infection and a high prevalence of chronic carriers.
According to numerous works (see for instance \cite{Coursaget, Edmunds, majid, McMahon, Peace} and the references therein) there is a rapid decline in the probability of developing the chronic carrier stage of the disease with respect to the age at infection. According Edmunds et al \cite{Ed93}, approximatively 90$\%$ of children of less than 5 years old will become carriers if infected while about 10$\%$ of adults will become carriers if infected.
Moreover this age-specific probability appears to be remarkably stable for a wide range of areas \cite{Ed93}.
Hence the mathematical modelling of HBV infection requires to take into account this 
strong age-specific differential susceptibility of humans.

In such highly endemic areas, these are two main routes for the infection of children: a vertical transmission from an infected mother to her infant, and an horizontal transmission through close contacts. 

Following WHO \cite{m1}, perinatal transmission is a much less important contributor to the carrier pool in Africa (we also refer to Assumption 4.1 in \cite{k1} for similar consideration for the Chinese case or \cite{abbas, majid} for the Pakistanese case). Moreover Anfumbom et al.\cite{anf} state that ``A low proportion of HBeAg among HBsAg-positive pregnant women with known HIV status could suggest low perinatal transmission of HBV for an endemic country like in Cameroon'' (see also  \cite{hw,kiire}). Hence, as in the work of Edmunds et al.\cite{Ed93,ccv}, we shall focus, in this work, on horizontal transmission taking into account age-specific susceptibility.

This differential susceptibility is a particularly important point for HBV infection. Let us recall that according to CDC\footnote{Centers for Disease Control and Prevention, USA: www.cdc.gov} (see also \cite{Ed93} and the references therein) about 90$\%$ of children will remain chronically infected with HBV while 90$\%$ of adults will develop acute infection and will completely recover from HBV infection.\\

Hence the main feature of the above described models consists in the age-specific susceptibility dependence through function $p_I(a)$  \cite{hw,kiire} (or $p_J(a)$). Edmunds et al. in \cite{Ed93} justifies the fact that ``the expected probability'' $p_J$ (of developing the chronic carrier state given the age $a$ at infection in an exponential function of age $a$) is decreasing. 
More specifically according Edmunds et al. in \cite{Ed93} function $q\equiv p_J$ takes the form $p_J(a)=p_1(a)=\kappa e^{-ra^s}$ for some suitable parameter set.
In order to take into account this age-specific susceptibility dependence we will use in this work a simplest prototypical shape curve of the form
\begin{equation}\label{prot1}
p_J(a)=\kappa e^{- ra },
\end{equation}
for some $\kappa \in [0,1]$ and $r>0$. Note that this specific form will only be used to check some assumptions in order to apply our main perturbation result. Based on the data of \cite{Ed93} we estimate the above parameters $\kappa$ and $r$ using the least squares method to obtain the following function (graphically closed to Edmunds's\cite[Fig. 1, page 199]{Ed93} function $p_1$ for $r=0.645$ and $s=0.455$):
\begin{equation}\label{function-q}
p_J(a)=p_2(a)=0.643\exp\left(-0.156a\right),\;a>0.
\end{equation}

The above models \eqref{eq1.1}-\eqref{eq1.2i} and \eqref{eq0}-\eqref{eq000} have been suggested in the work of Bonzi et al \cite{bonzi} using ordinary differential equations and discrete age structure for susceptible population. We refer to Zou et al \cite{k1} and the references therein for other classes of age structured system modelling HBV transmission.

Standard methodologies apply to provide the existence and uniqueness of mild solution for Systems \eqref{eq1}-\eqref{eq3} and \eqref{eq0}-\eqref{eq000} (see for instance Theorem \ref{THEO1} in section \ref{sect1f}, also \cite{Magal2009,Magal-Ruan,Thieme1990b} and the references therein).

The aim of this manuscript is to investigate the asymptotic behaviour of an ``aggregated'' system concerning infective compartments (with $I(t):=\int_0^\infty i(t,x)dx$ and $J(t):=\int_0^\infty j(t,x)dx$) derived from systems \eqref{eq1.1}-\eqref{eq1.2i} and \eqref{eq0}-\eqref{eq000}. To do so let us recall that according to WHO \cite{who}, Edmunds et al \cite{EMN1996},  Bonzi et al \cite{bonzi}, Fall et al \cite{cari} and Wilson et al \cite{WNC1998, WNC2000}, chronic carriers (most of time asymptomatic) have a low infectious rate. As a consequence in most part of this work we assume that
\begin{equation}\label{ass}
0\leq \beta_j<<\beta_i.
\end{equation}

In the extreme case when $\beta_j\equiv 0$ systems \eqref{eq0} or \eqref{eq1.1} re-write as the so-called reduced system (omitting the recovered class).

We shall assume that the contacts between individuals are homogeneous among the different cohorts, so that functions $\beta_i$ and $\beta_j$ are constant
\begin{equation*}
\beta_i\equiv \beta_I>0\text{ and }\beta_j\equiv \beta_J\geq 0.
\end{equation*}

 We will then study the general semiflow generating $(s,I,J)$ as a perturbation on $\beta_J$ around 0, of the reduced systems.

In addition for each function $p\in L^\infty(0,\infty)$ we set $p^*\in \left(L^1(0,\infty)\right)'$ the dual form associated to $p$ and defined by 
\begin{equation*}
p^*[\varphi]=\int_0^\infty p(a)\varphi(a)da,\;\;\forall \varphi\in L^1(0,\infty).
\end{equation*}

Our technical main \textbf{``aggregated'' system} \eqref{eq1}-\eqref{eq3} is obtained from both models \eqref{eq1.1}-\eqref{eq1.2i} and \eqref{eq0}-\eqref{eq000} by setting $
I(t):=\int_0^\infty i(t,x)dx$ and $J(t):=\int_0^\infty j(t,x)dx$, such that $(s,I,J)$ satisfies the following closed system of equations:
\begin{equation}\label{eq1}
\begin{cases}
\left(\partial_t+\partial_a+\mu\right)s(t,a)=-\lambda(t)s(t,a),\;\;t>0,\;a>0,\\
s(t,0)=\Lambda,\\
I'(t)=\lambda(t)p_I^*\left[s(t,.)\right]-\nu_I I(t),\\
J'(t)=\lambda(t)p_J^*[s(t,.)]-\nu_J J(t),\;\;t>0,
\end{cases}
\end{equation}
wherein we have set
\begin{equation}\label{eq2}
\lambda(t)=\beta_I I(t)+\beta_J J(t).
\end{equation}

Next by setting $I_0:=\int_0^\infty i_0(x)dx$ and $J_0:=\int_0^\infty j_0(x)dx$ one supplements \eqref{eq1} together with the initial data: 
\begin{equation}\label{eq3}
s(0,.)=s_0(.)\in L^1_+(0,\infty),\;\;I(0)=I_0\geq 0,\;\;J(0)=J_0\geq 0.
\end{equation}
We investigate the long time behaviour of \eqref{eq1}-\eqref{eq3} when $\beta_J$ is small using some knowledge about the dynamics of the auxiliary reduced systems with $\beta_J\equiv 0$.
The dynamics for the chronic carriers is fully known from the dynamics of the reduced system because the $(s,I)-$variables decouple from the $J-$variable.

In order to understand the dynamics of the above problem \eqref{eq1}-\eqref{eq3} with $0<\beta_J<<1$, we will make use  of the abstract results of Magal in \cite{jdemagal} to derive a result about global asymptotic stability of the endemic equilibrium for small values of $\beta_J$. These abstract results of Magal\cite{jdemagal} developed in the persistence context, the global perturbation of stable equilibrium theory initiated by Smith and Waltman\cite{hl}.

The work is organized as follows. In Section 2, we deal with the well-posedness of the system, derive preliminary results that will be useful to study the long term behaviour of the model and we study the global asymptotic stability of the disease free equilibrium in the case when the basic reproduction number satisfies $R_0\leq 1$ (see \eqref{Def-R0} for the definition of this threshold parameter). 
In Section 3 we state and prove our main perturbation result. It is concerned with the long time dynamics of \eqref{eq1}-\eqref{eq3} with $\beta_J<<1$ small enough.
It is stated under rather general assumptions without using the specific form of function $p_J$ while the proof is based on the perturbation arguments derived by Magal in \cite{jdemagal}. 
Finally Section 4 is devoted to an application of this perturbation result to our prototypical function \eqref{prot1} and concluding remarks end the manuscript.

\section{Preliminary results}\label{sect1f}

The aim of this section is to investigate basic properties of System \eqref{eq1}-\eqref{eq3}.
This analysis will be related to the following assumptions:
\begin{assumption}\label{ASS1}
We assume that $\Lambda>0$, $\mu>0$, $\beta_I>0$, $\beta_J\geq 0$, $\nu_I>0$ and $\nu_J>0$ are given parameters.
\end{assumption}
Our second assumption is related to function $p_I$.
\begin{assumption}\label{ASS2}
We assume that $p_I\in L^\infty_+(0,\infty)$ is a given measurable and bounded function such that
\begin{equation*}
0\leq p_I(a)\leq 1,\;for\;\;a.e.\;\;a\geq 0.
\end{equation*}
\end{assumption}
In order to deal with \eqref{eq1}-\eqref{eq3}, let us introduce the Banach spaces
\begin{equation*}
X=L^1(0,\infty)\times\R\times\R\times\R,\;\;X_0=L^1(0,\infty)\times\{0\}\times\R\times\R
\end{equation*}
endowed with the usual product norm, as well as its positive cone $X_+$ defined by
\begin{equation*}
X_+=L^1_+(0,\infty)\times [0,\infty)\times [0,\infty)\times [0,\infty)\text{ and }X_{0+}=X_0\cap X_+.
\end{equation*}
Consider the linear operator $A:D(A)\subset X\to X$ defined by 
\begin{equation*}
D(A)=W^{1,1}(0,\infty)\times\{0\}\times \R^2 \text{ and }A\begin{pmatrix} \varphi\\ 0\\\alpha_I\\\alpha_J\end{pmatrix}= 
\begin{pmatrix}
-\varphi'-\mu \varphi\\ -\varphi(0)\\ -\nu_I \alpha_I\\ -\nu_J \alpha_J
\end{pmatrix}.
\end{equation*}
Consider also the non-linear map of the class $C^\infty$, $F:X_{0+}\to X$ defined by
\begin{equation*}
F\begin{pmatrix} \varphi\\ 0\\\alpha_I\\\alpha_J\end{pmatrix}= 
\begin{pmatrix}
-(\beta_I \alpha_I+\beta_J\alpha_J)\varphi\\ \Lambda\\ (\beta_I \alpha_I+\beta_J\alpha_J)p_I^*[\varphi]\\(\beta_I \alpha_I+\beta_J\alpha_J)p_J^*[\varphi]
\end{pmatrix}.
\end{equation*}
Now identifying $\left(s(t,.),I(t),J(t)\right)$ together with $u(t)=\left(s(t,.),0,I(t),J(t)\right)^T$, System \eqref{eq1}-\eqref{eq3} re-writes as the following abstract Cauchy problem
\begin{equation}\label{Cauchy-pb}
\begin{split}
u'(t)=Au(t)+F\left(u(t)\right),\;\;t>0,\\
u(0)=\left(s_0(.),0,I_0,J_0\right)^T\in X_{0+}.
\end{split}
\end{equation}
Together with these notations, our first results are collected in the following theorem:
\begin{theorem}\label{THEO1}
Let Assumptions \ref{ASS1} and \ref{ASS2} be satisfied.
Then System \eqref{eq1}-\eqref{eq3} generates a strongly continuous positive semiflow $\{U(t)\}_{t\geq 0}$ on $X_{0+}$. This means that for each $x\in X_{0+}$, the continuous map $u:t\to u(t):=U(t)x$ defined from $[0,\infty)$ into $X_{0+}$ is a weak solution of \eqref{Cauchy-pb}, that is
\begin{equation*}
\begin{split}
&\int_0^t u(s)ds\in D(A),\;\;\forall t\geq 0,\\
&u(t)=x+A\int_0^t u(s)ds+\int_0^t F\left(u(s)\right)ds,\;\;\forall t\geq 0.
\end{split}
\end{equation*}
The semiflow $\{U(t)\}_{t\geq 0}$ satisfies the following properties
\begin{itemize}
\item [(i)] if $x=\left(s_0,0,I_0,J_0\right)^T\in X_{0+}$, then denoting $U(t)x=\left(s(t,.),0,I(t),J(t)\right)^T$, it satisfies the following Volterra integral formulation:
\begin{equation}\label{eq-Vol1}
s(t,a)=\begin{cases} s_0(a-t)e^{-\int_0^t [\mu+\lambda(z)]dz}\text{ if $a\geq t\geq 0$}\\
\Lambda e^{-\int_{t-a}^{t} [\mu+\lambda(z)]dz}\text{  \;\;\;\;\;\;\; if $0\leq a<t$},
\end{cases}
\end{equation}
while functions $t\to I(t)$ and $t\to J(t)$ are of class $C^1$ on $[0,\infty)$ and satisfy
\begin{equation}\label{eq-Vol2}
\begin{cases}
I(0)=I_0,\;\;J(0)=J_0,\;\;\lambda(t)=\beta_I I(t)+\beta_J J(t),\;\;\forall t\geq 0,\\
I'(t)=\lambda (t)p_I^*[s(t,.)]-\nu_I I(t),\\
J'(t)=\lambda(t)p_J^*[s(t,.)]-\nu_J J(t),\;\;t>0.
\end{cases}
\end{equation}
\item [(ii)] It satisfies the following bounded-dissipative estimates for each $x\in X_{0+}$ and each $t\geq 0$:
\begin{equation*}
\|U(t)x\|_{X_0}\leq \frac{\Lambda}{\nu}\left(1-e^{-\nu t}\right)+\|x\|_X e^{-\nu t},
\end{equation*} 
wherein we have set $\nu=\min\left\{\mu,\nu_{I},\nu_{J}\right\}$.
\item [(iii)] The semiflow $\{U(t)\}_{t\geq 0}$ is asymptotically smooth on $X_0$.
\end{itemize}
\end{theorem}

\begin{proof}
The proof of $(i)$ is rather standard. Indeed it is easy to
check that operator $A$ satisfies the Hille-Yosida property. Then
standard methodologies apply to provide the existence and
uniqueness of mild solution for System \eqref{eq1}-\eqref{eq3}.
(see for instance \cite{Magal2009,
Magal-Ruan,Thieme1990b} and the references therein).\\
The proof of $(ii)$ is immediate from the integration of the equations.\\
We will now focus $(iii)$. In order to prove this result we will make use of results derive by Sell and You in  \cite{Sell-2002}. More precisely, we will show that for each bounded set $B\subset X_{0+}$, the semiflow $\{U(t)\}_{t\geq 0}$ is asymptotically compact on $B$. 
Let $B\subset X_{0+}$ be a given bounded set.
Let us consider a sequence $\left(s_0^j,0,I_0^j, J_0^j\right)_{j\geq 0}\subset B$ of initial data and let us denote by $$\left(s^j(t,.),0,I^j(t),J^j(t)\right)^{T}=U(t)\left(s_0^j,0,I_0^j, J_0^j\right)^T$$ the solution semiflow.
Let $\{t_j\}_{j\geq 0}$ be a sequence tending to $\infty$. We aim to show that the sequence $\left(s^j(t_j,.),0,I^j(t_j),J^j(t_j)\right)^{T}$ is relatively compact in $X_0$. To that aim, let us first notice that due to estimate $(ii)$ in Theorem \ref{THEO1}, 
the sequence of functions $J^j(t_j+.)$ and $I^j(t_j+.)$ are uniformly bounded as well as their time derivatives.
Using the Arzela Ascoli theorem, possibly along a subsequence, one may assume that 
$(J^j(t_j+.),I^j(t_j+.))$ converges locally and uniformly to some bounded, continuous and positive functions $t\in\R\to (J(t),I(t))$. 
Then setting for each $j\geq 0$, $\lambda_j(t)=\beta_I I^j(t)+\beta_J J^j(t)$ and
$\lambda:\R\to [0,\infty)$ defined by $\lambda(t)=\beta_I I(t)+\beta_J J(t)$, one obtains using Lebesgue convergence theorem that as $j\to\infty$
\begin{equation*}
\int_{t_j}^\infty s_0^j (a-t_j) da e^{-\int_0^{t_j} [\mu+\lambda_j(s)]ds}\leq \int_0^\infty s_0^j(a)da e^{-\mu t_j}\to 0.
\end{equation*}
while as $j\to\infty$
\begin{equation*}
e^{-\int_{t_j-a}^{t_j} [\mu+\lambda_j(s)]ds}=e^{-\int_{-a}^{0} [\mu+\lambda_j(t_j+s)]ds}\to e^{-\int_{-a}^{0} [\mu+\lambda(s)]ds}\;\;for\; \;a.e.\;\,a\geq 0. 
\end{equation*}
As a consequence one obtains that
\begin{equation*}
s_j(t_j,a)=s_0^j(a-t_j)e^{-\int_0^{t_j} [\mu+\lambda_j(s)]ds}{\bf 1}_{\{a\geq t_j\}}+\Lambda e^{-\int_{t_j-a}^{t_j} [\mu+\lambda_j(s)]ds}{\bf 1}_{\{a<t_j\}}
\end{equation*}
and
\begin{equation*}
\lim_{j\to\infty}s_j(t_j,.)=\Lambda e^{-\int_{-a}^{0} [\mu+\lambda(s)]ds}\text{ in $L^1(0,\infty)$}.
\end{equation*}
The result follows.
\end{proof}

From Theorem \ref{THEO1} one deduces using the results of Hale \cite{Hale-1989} (see also Smith and Thieme \cite{Smith-Thieme}, Magal and Zhao \cite{Magal-Zhao} and the references therein), the following results:
\begin{proposition}\label{PROP3}
Let Assumption \ref{ASS1}-\ref{ASS2} be satisfied. The semiflow $\{U(t)\}_{t\geq 0}$ provided by Theorem \ref{THEO1} has a non-empty compact global attractor $\mathcal A\subset X_{0+}$. This means that $\mathcal A$ is compact, invariant and attracts all bounded set $B\subset X_{0+}$, in the sense that for each $B\subset X_{0+}$ bounded subset, one has $d\left(U(t)B,\mathcal A\right)\to 0$ as $t\to\infty$ where $d(B,A)$ denotes the semi distance from $B$ to $A$ defined by
\begin{equation*}
d(B,A)=\sup_{y\in B}\inf_{x\in A} \|y-x\|_X.
\end{equation*} 
\end{proposition}
Due to the above result the description of the long time behaviour of \eqref{eq1}-\eqref{eq3} relies on the suitable descriptions of the global attractor. 
Now we need to derive some basic properties of the entire solutions of \eqref{eq1}-\eqref{eq3}.
Before doing so, let us recall that the entire solutions are smooth, in the sense that $s\in C^\infty(\R\times [0,\infty))$ while $(I,J)\in C^\infty(\R)^2$.
With this in mind, we can state our first lemma.
\begin{lemma}\label{LE-upper}
Assume that the assumptions of Proposition \ref{PROP3} are satisfied.
Let $t\in \R\mapsto \left(s(t,.),I(t),J(t)\right)^T$ be a given entire solution of \eqref{eq1}-\eqref{eq3}, then the following property holds true:
\begin{equation}\label{basic-upper}
s(t,a)\leq \Lambda e^{-\mu a}:= s_F(a),\;\forall t\in\R,\;\forall \;a\geq 0.
\end{equation}
One furthermore has:
\begin{equation}\label{aqz1}
\left(\exists t_0\in \R\;\;\beta_I I(t_0)+\beta_J J(t_0)=0\right)\;\Rightarrow\;\left(I(t)\equiv J(t)\equiv 0,\forall t\in\R\right)
\end{equation}
and,
\begin{equation}\label{aqz2}
\left(\exists t_0\in \R,\;a_0>0\;\;s(t_0,a_0)\equiv \Lambda e^{-\mu a_0}\right)\;\Rightarrow\;\begin{cases}I(t)\equiv J(t)\equiv 0,\forall t\in\R,\\
s(t,.)\equiv s_F(.),\;t\in\R.\end{cases}
\end{equation}
\end{lemma}

\begin{proof}
In order to prove \eqref{basic-upper} we make use of the integral Volterra formulation of the solution that reads as follow: for each $t\in\R$ and $a>0$
\begin{equation}\label{formule}
s(t,a)=\Lambda e^{-\mu a-\int_{t-a}^{t}\lambda(l)dl}.
\end{equation}
Hence since $\lambda \geq 0$, \eqref{basic-upper} follows.

Let us now prove \eqref{aqz1}. With this in mind, denote by $\overline\nu:=\max\left(\nu_I,\nu_J\right)$ and $\underline\nu:=\min\left(\nu_I,\nu_J\right)$. Then due to \eqref{eq1} one has:
\begin{equation*}
\lambda(t)\left[g(t)-\overline{\nu}\right]\leq \lambda'(t)\leq \lambda(t)\left[g(t)-\underline{\nu}\right],\;\;\forall t\in\R,
\end{equation*}
where we have set $g(t):=\int_0^\infty \left[\beta_I p_I(a)+\beta_J p_J(a)\right]s(t,a)da$.
Then if there exists $t_0\in\R$ such that $\lambda(t_0)=0$ then due to the above differential inequality one obtains $\lambda(t)\equiv 0$.
As a consequence, since $\beta_I>0$, on the one hand if $\beta_J>0$ then the result follows.
On the other hand, if $\beta_J=0$ then one gets from the above that $I(t)\equiv 0$ and $J$ becomes a uniformly bounded entire solution of $J'(t)=-\nu_J J(t)$ for all $t\in\R$. Hence $J(t)\equiv 0$.

Finally to complete the proof of \eqref{aqz2}, Note that if there exist $t_0\in\R$ and $a_0>0$ with $s(t_0,a_0)=s_F(a_0)$ then using \eqref{formule}, there exists $t_1\in \left(t_0-a_0,t_0\right)$ such that $\lambda(t_1)=0$. Thus due to \eqref{aqz1}, $\lambda(t)\equiv 0$ while \eqref{formule} ensures $s(t,.)\equiv s_F(.)$. The result follows. 
\end{proof}

In order to give some information on the long term behaviour of system \eqref{eq1}-\eqref{eq3}, we introduce the threshold parameter $R_0$ defined by
\begin{equation}\label{Def-R0}
R_0:=\Lambda \left[\frac{\beta_I}{\nu_I}p_I^*\left[e^{-\mu.}\right]+\frac{\beta_J}{\nu_J}p_J^*\left[e^{-\mu.}\right]\right].
\end{equation}
Using this definition, straightforward computations lead us to the following result:
\begin{lemma}\label{LE4}
Let Assumption \ref{ASS1}-\ref{ASS2} be satisfied. Then the following holds true:
\begin{itemize}
\item [(i)] If $R_0\leq 1$, then System \eqref{eq1}-\eqref{eq3} has a unique stationary state $x_F=\left(s_F,0,0,0\right)^T\in X_+$ where $s_F(a)=\Lambda e^{-\mu a}$.\\
\item [(ii)]If $R_0>1$, then system \eqref{eq1}-\eqref{eq3} has two stationary states $x_F\in X_+$ and $x_E=\left(s_E(.),0,I_E, E_E\right)^T$ with
\begin{equation*}
s_E(a)=s_F(a) e^{-\lambda_E a},\;\;I_E=\frac{\lambda_E}{\nu_I}p_I^*\left[s_E\right],\;\;J_E=\frac{\lambda_E}{\nu_J}p_J^*\left[s_E\right],
\end{equation*}
and where $\lambda_E>0$ is the unique solution of the equation
\begin{equation*}
1=\Lambda \sum_{k\in\{I,J\}}\frac{\beta_k}{\nu_k}p_k^*\left[e^{-\left(\mu+\lambda_E\right).}\right].
\end{equation*}
\end{itemize}
\end{lemma}
Together with this threshold definition, we are now able to study the global dynamics of \eqref{eq1}-\eqref{eq3} when $R_0\leq 1$.
Our result is stated in the next proposition:
\begin{proposition}[The case $R_0\leq 1$: disease extinction]\label{PROP4}
Let Assumption \ref{ASS1}-\ref{ASS2} be satisfied. Assume that $R_0\leq 1$, then the global attractor $\mathcal A$ defined in Proposition \ref{PROP3} satisfies $\mathcal A=\left\{x_F\right\}$.
\end{proposition}

\begin{proof}
Let us consider for $k\in \{I,J\}$ the quantities $\Gamma_{k}>0$ defined by $\Gamma_{k} = \frac{\beta_{k}}{\nu_{k}}$.
Using these notations, note that one can re-write
\begin{equation*}
R_{0}=\sum_{k\in\{I,J\}}\Gamma_{k}p_k^*\left[s_F\right].
\end{equation*}
Next consider the map $L:X_+\to \R$ defined by
\begin{equation*}
L\left[\left(s, 0, I, J\right)^T\right]=\Gamma_{I}I+\Gamma_{J}J.
\end{equation*}
Let $x\in \mathcal A$ be given and let us denote by $\left\{V(t)=\left(s(t,.),0,I(t),J(t)\right)^T\right\}_{t\in\R}\subset \mathcal A$ an entire solution of \eqref{eq1} with $V(0)=x$.
Then straightforward computations yield the following: for each $t\in\R$ and $\tau\geq 0$:
\begin{equation}\label{lyapunov1}
L\left[V(t+\tau)\right]-L\left[V(t-\tau)\right]=\int_{t-\tau}^{t+\tau} \left[ \sum_{k\in\{I,J\}} \Gamma_{k}p_k^*\left[s(l,.)\right]-1\right]  \left[ \beta_{I}I(l)+\beta_{J}%
J(l)\right]dl
\end{equation}
Recalling that $R_0\leq 1$ and using \eqref{basic-upper}, one obtains that the map $t\mapsto L\left[V(t)\right]$ is bounded and decreasing on $\R$.\\
Let $\{t_n\}_{n\geq 0}$ be a decreasing sequence tending to $-\infty$ as $n\to\infty$. Consider the sequence of maps $\left\{V_n=V(.+t_n)\right\}_{n\geq 0}$. Up to a subsequence one may assume that $V_n$ converges towards some function $\widehat{V}=\left(\widehat s, 0,\widehat I,\widehat J\right)^T$ for the topology of $C_{loc}\left(\R,X\right)$. Note that $\left\{\widehat{V}(t)\right\}_{t\in\R}$ is also an entire solution of \eqref{eq1}-\eqref{eq3}.
Using \eqref{lyapunov1}, one obtains that $\widehat V$ satisfies:
\begin{equation}\label{lyapunov2}
\left[  \sum_{k\in\{I,J\}}\Gamma_kp_k^*\left[\widehat s(t,.)\right]-1\right]  \left[  \beta_{I}\widehat{I}(t)+\beta_{J}%
\widehat{J}(t)\right]\equiv 0.
\end{equation}
Furthermore since $\{t_n\}$ is decreasing and the map $t\mapsto L\left[V(t)\right]$ is decreasing, this leads us to
\begin{equation}\label{lyapunov3}
0\leq L\left[V(t)\right]\leq L\left[\widehat{V}(0)\right],\;\;\forall t\in\R.
\end{equation}
Now if $R_0<1$, we infer from \eqref{lyapunov2} that $\widehat{I}(t)\equiv 0$ and $\widehat{J}(t)\equiv 0$. 
From \eqref{lyapunov3}, one gets that $I(t)\equiv J(t)\equiv 0$ and since $V$ is an entire solution of \eqref{eq1}-\eqref{eq3}, one obtains that $V(t)\equiv x_F$ and the result follows.\\
If we assume that $R_0=1$, then, by using Lemma \ref{LE-upper}, \eqref{aqz1} and \eqref{aqz2}, \eqref{lyapunov2} implies that either $\widehat{s}(t,.)\equiv s_F$ or $I(t)\equiv J(t)\equiv 0$. The second condition can be handled similarly to the argument for $R_0<1$. In the first case, namely if $\widehat{s}(t,.)\equiv s_F$, one gets from the $s-$equation that
\begin{equation*}
\left(\beta_I \widehat{I}(t)+\beta_J \widehat{J}(t)\right)s_F(a)\equiv 0,
\end{equation*}
so that $\beta_I \widehat{I}(t)+\beta_J \widehat{J}(t)\equiv 0$ and the result follows.
\end{proof}

\section{The case $R_0>1$: a perturbation result}
In this section we will study global dynamical properties of \eqref{eq1}-\eqref{eq3} in the framework of the biological assumption  
$\beta_J<<1$ being a small parameter. For notational simplicity we replace $\beta_J$ by $\varepsilon$.
Moreover we will explicitly write down the dependence of several quantities with respect to $\varepsilon\geq 0$.
For instance, recalling \eqref{Def-R0}, we write $R_0=R_0[\varepsilon]$, $U(t)=U_\varepsilon(t)$ the semiflow provided by Theorem \ref{THEO1}, $x_E=x_E^\varepsilon$ the endemic equilibrium defined in Lemma \ref{LE4} and so on.

Before stating our main result, let us state our assumptions that will be discussed in the next section using the specific form of function $p_J$ in \eqref{prot1}:
\begin{assumption}\label{ASS3}
We assume that the following properties hold true:
\begin{itemize}
\item [(i)] $R_0[0]>1$ and,
\item[(ii)]  The semiflow $\left\{U_0(t)\right\}_{t\geq 0}$ satisfies:
\begin{equation*}
\lim_{t\to\infty} U_0(t)x=x_E^0,\;\;\forall x\in M_{*}:=\left\{\left(s_0,0, I_0, J_0\right)^T\in X_+:\;\;I_0>0\right\}.
\end{equation*}
\item [(iii)] The map $\Delta:\Omega\to \mathbb C$ defined by
\begin{equation*}
\begin{split}
&\Omega=\left\{z\in\mathbb C:\;{\rm Re}\;(z)>-\mu\right\}\\
&\Delta(\lambda)=\lambda+I_{E}\beta_{I}^{2}p_I^*\left[s_{E}(.)\left[
\frac{1-e^{-\lambda .}}{\lambda}\right]\right],
\end{split}
\end{equation*}
does not have any root with ${\rm Re}(\lambda)\geq 0$.
\end{itemize}
\end{assumption}

Then the following result holds true:
\begin{theorem}\label{THEO5} 
Let Assumption \ref{ASS1}, \ref{ASS2} and \ref{ASS3} be satisfied. Then there exists $\delta>0$ such that for each $\varepsilon\in (0,\delta)$, the semiflow $\{U_\varepsilon(t)\}_{t\geq 0}$ satisfies
\begin{equation*}
\lim_{t\to\infty} U_\varepsilon(t)x=x_E^\varepsilon,
\end{equation*}
for each $x=\left(s_0,0, I_0, J_0\right)^T\in X_{0+}$ such that $I_0+J_0>0$.
\end{theorem}

The proof of this result relies on the application of Theorem 1.2 derived by Magal in \cite{jdemagal} (see also Smith and Waltman \cite{hl} for an earlier result).
In order to use the above mentioned result, let us check the required set of assumptions. 
They are collected in the following lemma:
\begin{lemma}\label{LE-crucial}
Let Assumption \ref{ASS1}, \ref{ASS2} and \ref{ASS3} $(i)$ be satisfied. Consider the map $\rho:[0,\infty)\times X_{0+}\to [0,\infty)$ defined by
\begin{equation*}
\rho(\varepsilon,x)\equiv\rho_\varepsilon (x)=\beta_I I+\varepsilon J,\;\;\forall x=\left(s,0,I,J\right)^T\in X_{0+},
\end{equation*}
and let us set for each $\varepsilon\geq 0$:
\begin{equation}\label{defM00}
M_0^\varepsilon=\left\{x\in X_{0+}:\;\rho_\varepsilon(x)>0\right\},\;\;\partial M_0^\varepsilon=\left\{x\in X_{0+}:\;\rho_\varepsilon(x)=0\right\}
\end{equation}
Then the following holds true:
\begin{itemize}
\item [(i)] For each $\varepsilon\geq 0$, $M_0^\varepsilon$ and $\partial M_0^\varepsilon$ are both positively invariant under the semiflow $\{U_\varepsilon(t)\}_{t\geq 0}$,
\item [(ii)] For each $\varepsilon\geq 0$, $U_\varepsilon$ has a global attractor $\mathcal A_\varepsilon$ in $X_+$ and the family $\left\{\mathcal A_\varepsilon\right\}_{\varepsilon \geq 0}$ is upper semi-continuous at $\varepsilon=0$.
\item [(iii)] For each $\delta_0>0$, there exists $\eta>0$ such that for each $\varepsilon\in [0,\delta_0]$ and each $x\in M_0^\varepsilon$, one has
\begin{equation*}
\liminf_{t\to\infty}\rho_\varepsilon\left(U_\varepsilon(t)x\right)\geq \eta.
\end{equation*}
\end{itemize}
\end{lemma}

Before proving this result let us first complete the proof of Theorem \ref{THEO5}.
As already mentioned the proof of this result is a direct application of Theorem 1.2 in Magal \cite{jdemagal}.
To that aim, let us introduce for each $\varepsilon>0$ the $C^\infty-$map $F_\varepsilon:X_{0+}\to X$ defined by $F_\varepsilon=F_0+\varepsilon G$ where $F_0:X_{0+}\to X$ and $G:X_{0+}\to X$ are defined by
\begin{equation*}
F_0\begin{pmatrix} \varphi\\ 0\\\alpha_I\\\alpha_J\end{pmatrix}= 
\begin{pmatrix}
-\beta_I \alpha_I\varphi\\ \Lambda\\ \beta_I \alpha_Ip_I^*\left[\varphi\right]\\\beta_I \alpha_Ip_J^*\left[\varphi\right]
\end{pmatrix},\;\;G\begin{pmatrix} \varphi\\ 0\\\alpha_I\\\alpha_J\end{pmatrix}= 
\begin{pmatrix}
-\alpha_J \varphi\\ 0\\ \alpha_Jp_I^*\left[\varphi\right]\\ \alpha_Jp_J^*\left[\varphi\right]
\end{pmatrix}.
\end{equation*}

Let us introduce the linear semiflow $\{L(t)\}_{t\geq 0}\subset \mathcal L\left(X_0\right)$ defined for each $x\in X_0$ by $L(t)x=u(t)$ where $u(t)$ is the mild solution of the following linear Cauchy problem:
\begin{equation}\label{eq-linear}
\frac{du(t)}{dt}=Au(t)+DF_0\left(x_E^0\right)u(t),\;t>0\text{ and }u(0)=x.
\end{equation}
We first claim that:
\begin{lemma}\label{claim1}
For each $t>0$ the spectral radius of $L(t)$ denoted by $r(L(t))$ satisfies:
\begin{equation*}
r(L(t))<1.
\end{equation*}
\end{lemma}
The proof of this lemma is postponed.\\

Let us denote for any subset $A\subset X_{0+}$ and map $f:A\to X$:
\begin{equation}\label{DEF-Lip}
\|f\|_{{\rm Lip},\;A}:=\sup_{x,y\in A,\;x\neq y}\frac{\|f(x)-f(y)\|}{\|x-y\|}.
\end{equation}
Using this notation we claim that:
\begin{lemma}\label{claim2}
For any given and fixed $t_0>0$ one has:
\begin{equation}\label{estimation}
\lim_{\delta\to 0^+} \sup_{\varepsilon \in [0,\delta]}\left\|U_\varepsilon(t_0)-L(t_0)\right\|_{{\rm Lip},\overline{B}_{X_{0}}\left(x_E^0,\delta\right)\cap X_{0+}}=0.
\end{equation} 
\end{lemma}
Before proving this lemma, let us notice that Theorem \ref{THEO5} directly follows from the properties stated in Lemma \ref{LE-crucial}, Lemma \ref{claim1} and Lemma \ref{claim2}

The rest of this section is devoted to the proof of Lemma \ref{LE-crucial}, Lemma \ref{claim1} and Lemma \ref{claim2}. 

\subsection{Proof of Lemma \ref{LE-crucial}}
This section is devoted to the proof of Lemma \ref{LE-crucial}.\\
Let us first notice that $(i)$ is straightforward.\\
In order to prove $(ii)$ we will make use of Proposition 2.9 in \cite{jdemagal}.
Recall that due to Proposition \ref{PROP3}, for each $\varepsilon\geq 0$ the semiflow $\{U_\varepsilon(t)\}_{t\geq 0}$ is bounded, dissipative, asymptotically smooth and has a global attractor denoted by $\mathcal A_\varepsilon\subset X_{0+}$ that attracts all bounded subsets. According to Proposition 2.9 in \cite{jdemagal}, in order to prove that the family $\{\mathcal A_\varepsilon\}$ is upper semi-continuous at $\varepsilon=0$ it is sufficient to show that
\begin{equation}\label{conv}
U_\varepsilon(t)x\to U_0(t)x,\text{ as }\varepsilon\to 0^+,
\end{equation}
uniformly with respect to $(t,x)$ on bounded sets.
This is an application of Gronwall's Lemma.
Let $\tau>0$ and $\kappa>0$ be given.
Set $B=\overline{B}_{X_0}(0,\kappa)\cap X_{0+}$ and note that due to Theorem \ref{THEO5} $(ii)$ one has for all $\varepsilon\geq 0$ and $t\geq 0$:
\begin{equation*}
U_\varepsilon(t)B\subset \tilde B\text{ with }\tilde B=\overline{B}_{X_0}\left(0,\frac{\Lambda}{\mu}+\kappa\right)\cap X_{0+}.
\end{equation*}
Denote by $A_0:D(A_0)\subset \overline{D(A)} \rightarrow \overline{D(A)}$ the part of $A$ in $X_0=\overline{D(A)}$, which is defined by
\begin{equation}\label{part}
A_0 x=Ax,\; \forall x\in D(A_0)=\left\{x\in D(A) \; : \; Ax\in \overline{D(A)}\right\}.
\end{equation}
Recall that $A_0$ is densely defined and generates a $C^0$-semigroup denoted by $\left\{T_{A_0}(t)\right\}_{t\geq 0}$ in $X_0$. Here and in the sequel, we will denote by $R(\lambda;A)=(\lambda I-A)^{-1}$ the resolvent of $A$ for each value $\lambda\in \rho(A)$, the resolvent set of $A$.

Let $x\in B$ be given. Then for each $\varepsilon\geq 0$ and $t\geq 0$ one has:
\begin{equation*}
U_\varepsilon(t)x=T_{A_0}(t)x+\lim_{\lambda\to \infty}\int_0^t T_{A_0}(t-s)\lambda R(\lambda;A)F_\varepsilon\left(U_\varepsilon(s)\right)ds.
\end{equation*}
We refer to Magal and Ruan \cite{Magal-Ruan} for details on the above constant variation formula.
Hence, since $A$ satisfies the Hille-Yosida property, there exists some constant $M>0$ such that for all $t\in [0,\tau]$:
\begin{equation*}
\|U_\varepsilon(t)x-U_0(t)x\|\leq \varepsilon M\sup_{y\in \tilde B}\|G(y)\|+M\|F_0\|_{{\rm Lip},\;\tilde B}\int_0^t \|U_\varepsilon(s)x-U_0(s)x\|ds.
\end{equation*}
Then by using Gronwall's inequality, one completes the proof of \eqref{conv} and $(ii)$ follows.

Next it remains to prove $(iii)$. 
To that aim we will prove uniform weak persistence and will convert later into uniform strong persistence\cite{Thieme2}. We include parameter $\varepsilon$ into the state space.
We consider the Banach space $Y=X\times \R$ and define:
\begin{equation*}
Y_0=X_0\times \R,\;Y_+=X_+\times \R^+\text{ and }Y_{0+}=Y_0\cap Y_+.
\end{equation*}
Then we consider the following problem:
\begin{equation}\label{abstract1}
\begin{cases}
u'(t)=Au(t)+F_0(u(t))+\varepsilon G\left(u(t)\right),\;t>0,\\
\varepsilon'(t)=0,
\end{cases}
\end{equation}
supplemented together with initial data $\left(u(0),\varepsilon(0)\right)^T\in Y_{0+}$.
Let us consider the semiflow $\{V(t)\}_{t\geq 0}$ on $Y_{0+}$ generated by \eqref{abstract1}.
It is defined by
\begin{equation*}
V(t)\begin{pmatrix}x\\ \varepsilon\end{pmatrix}=\begin{pmatrix} U_\varepsilon(t)x\\ \varepsilon\end{pmatrix},\;\;\forall \begin{pmatrix}x\\ \varepsilon\end{pmatrix}\in Y_{0+}.
\end{equation*} 
Let $\delta_0>0$ be given. Let us also consider the complete metric space 
$M:=X_{0+}\times\lbrack0,\delta_0]$, the continuous function $\widehat{\rho
}:M\rightarrow\lbrack0,\infty)$ defined by
\begin{equation}
\widehat{\rho}(y)=I+\varepsilon J,\;\;y=\left(s,0,I,J,\varepsilon\right)\in M,
\end{equation}
as well as the sets
\begin{equation}
M_0=\{y\in M:\widehat{\rho}(y)>0\}\text{ and }\partial M_0=M\setminus M_0.
\end{equation}
Then the following result holds true: 
\begin{lemma}\label{weak}
Under the above assumptions and using the above notations the semiflow $\{V(t)\}_{t\geq 0}$ is $\widehat \rho$ weakly uniformly persistent, namely  
there exists $\eta>0$ such that
\begin{equation*}
\limsup_{t\rightarrow\infty}\;\widehat{\rho}\left(  V(t)y\right)  \geq\eta,\;\;\forall y\in M_0.
\end{equation*}
\end{lemma}

\begin{proof}
In order to prove this result, let us first notice that due to Theorem \ref{THEO1}, for each $y\in M_0$ there exists a time $t_y>0$ such that
\begin{equation*}
V(t)y\in M_0\cap \left[\overline{B}_{X_0}\left(0,\frac{\Lambda}{\mu}+1\right)\times [0,\delta_0]\right],\;\;\forall t\geq t_y.
\end{equation*}
As a consequence to prove Lemma \ref{weak}, it is sufficient to show that for each $r>0$ there exists $\eta_r>0$ such that
\begin{equation*}
\limsup_{t\rightarrow\infty}\;\widehat{\rho}\left(  V(t)y\right)  \geq\eta_r,\;\;\forall y\in M_0(r),
\end{equation*}
wherein we have set $M_0(r)=M_0\cap \left[\overline{B}_{X_0}\left(0,r\right)\times [0,\delta_0]\right]$.\\
Let $r>0$ be given.
Let us assume by contradiction that there exists a sequence $\{y_n=\left(x_n,\varepsilon_n\right)\}_{n\geq 0}\subset  M_0(r)$ and a sequence $\{t_n\}_{n\geq 0}$ tending to $\infty$ as $n\to\infty$ such that:
\begin{equation*}
\widehat{\rho}\left(V(t_n+t)y_n\right)\leq \frac{1}{n+1},\;\;\forall n\geq 0,\;t\geq 0. 
\end{equation*}
Up to a subsequence, one may assume that $\varepsilon_n\to \varepsilon_0\in [0,\delta_0]$ and possibly up to extraction of an other subsequence one easily checked for the system that:
\begin{equation}\label{convergence-V}
\lim_{n\to\infty} V(t_n+t)y_n=y_\infty:=\begin{pmatrix} x_F\\ \varepsilon_0\end{pmatrix},
\end{equation}
uniformly with respect to $t\geq 0$. \\
Next, as a consequence of \eqref{convergence-V}, one obtains that for each $\kappa>0$ small enough there exists $n=n(\kappa)\geq 0$ large enough such that for any $k\in\left\{I,J\right\}$:
\begin{equation}\label{esti-a}
p_k^*\left[s_F\right]-\kappa\leq p_k^*\left[s^n(t,.)\right]\leq p_k^*\left[s_F\right]+\kappa,\;\;\forall t\geq 0.
\end{equation}
Here we have set for each $t\geq 0$ and $n\geq 0$.
\begin{equation*}
\left(s^n(t,.),0,I^n(t),J^n(t),\varepsilon_n\right)^T=V(t_n+t)y_n.
\end{equation*}
Note now that for each $n\geq 0$ and each $t\geq 0$ one has:
\begin{equation}\label{esti-b}
\left(\beta_I I^n(t)+\varepsilon_n J^n(t)\right)\geq \beta_I I^n(t).
\end{equation}
Let $\kappa>0$ be given such that $p_I^*\left[s_F\right]-\kappa>0$. Let $n=n(\kappa)$ be defined by \eqref{esti-a}.
Then we infer from \eqref{esti-a}, \eqref{esti-b} and the $I-$equation that for the above defined  value of $n$, namely $n=n(\kappa)$, one has for each $t\geq 0$
\begin{equation*}
\left(I^n\right)'(t)\geq \beta_I I^n(t)\left[p_I^*\left[s_F\right]-\kappa\right]-\nu_I I^n(t),\\
\end{equation*}
This re-writes as
\begin{equation*}
I^n(t)\geq I^n(0)e^{\varrho t},\;\;\forall t\geq 0,
\end{equation*}
with 
\begin{equation*}
\varrho=\beta_I\left[\,p_I^*\left[s_F\right]-\kappa\right]-\nu_I=\nu_I\left(R_0[0]-1-\frac{\beta_I}{\nu_I}\kappa\right)>0,
\end{equation*}
as soon as $\kappa>0$ is small enough, namely $\frac{\nu_I}{\beta_I}\left(R_0[0]-1\right)>\kappa$.  Here recall that one has assumed $R_0[0]>1$ in Assumption \ref{ASS3}-(i). 
Since $I^n(0)>0$ this last sub-estimate contradicts the boundedness of $I$.
This completes the proof of Lemma \ref{weak}.
\end{proof}

Let us recall that the semiflow $\{V(t)\}_{t\geq 0}$ is bounded dissipative on $M$ and it is asymptotically smooth so that the semiflow has a compact attractor that attracts all bounded subset of $M$. (see Theorem 2.9 in \cite{Magal-Zhao}).
Then due to Proposition 3.2 in Magal and Zhao \cite{Magal-Zhao} or Theorem 2.3 of Thieme\cite{Thieme2}, Lemma \ref{weak} ensures that Lemma \ref{LE-crucial} $(iii)$ holds true.
This completes the proof of Lemma \ref{LE-crucial}.

\subsection{Proof of Lemma \ref{claim1}}

In order to prove that Lemma \ref{claim1} holds true we investigate some spectral properties of the Hille-Yosida operator $A+DF_0\left(x_E^0\right):D(A)\subset X\to X$.
Note that one has:
\begin{equation*}
DF_0\left(x_E^0\right)\begin{pmatrix}\varphi\\ 0\\ \alpha_I\\ \alpha_J\end{pmatrix}=\begin{pmatrix} -\beta_I s_E(.)\alpha_I-\beta_I I_E\varphi\\ 0\\ \alpha_I\beta_I p_I^*\left[s_E\right] +\beta_I I_E p_I^*\left[\varphi\right]\\
\alpha_I \beta_I p_J^*\left[s_E\right]+\beta_I I_E p_J^*\left[\varphi\right]
\end{pmatrix}.
\end{equation*}
Hence it is easy to check that (see \cite{Ducrot08})
\begin{equation*}
\omega_{0,ess}\left(\left(A+DF_0\left(x_E^0\right)\right)_{X_0}\right)\leq -\left[\mu+\beta_I I_E\right]\leq -\mu,
\end{equation*}
where $\omega_{0,ess}$ denotes the essential growth rate while the index $X_0$ corresponds to the part in $X_0$ of the linear operator.
As a consequence of this remark, if $\sigma\left(A+DF_0\left(x_E^0\right)\right)$ denotes the spectrum of the linear operator and recalling the definition of $\Omega$ in Assumption \ref{ASS3}, then one obtains using the results in \cite{Engel, Webb} that 
\begin{equation*}
\sigma\left(A+DF_0\left(x_E^0\right)\right)\cap\Omega
\text{ is only composed of point spectrum}.
\end{equation*}
As a consequence of the spectral mapping theorem, since the linear semigroup $\{L(t)\}_{t\geq 0}$ is generated by the linear operator $\left(A+DF_0\left(x_E^0\right)\right)_{0}$, the part of $A+DF_0\left(x_E^0\right)$ in $X_0$ (see definition in \eqref{part}), in order to prove Lemma \ref{claim1}, it is sufficient to prove that
\begin{equation}\label{eq-suff}
\sigma\left(A+DF_0\left(x_E^0\right)\right)\cap\Omega\subset \{\lambda\in\Omega:\;{\rm Re}\left( \lambda\right)<0\}.
\end{equation}
Let us notice that, since $DF_0\left(x_E^0\right)X_0\subset X_0$, one has $\left(A+DF_0\left(x_E^0\right)\right)_{0}=A_{0}+DF_0\left(x_E^0\right)_{|X_0}$.\\
Now it is easy to check that $0\not\in\sigma\left(A+DF_0\left(x_E^0\right)\right)$.
Next let $\lambda\in\Omega\setminus\{0\}$ and $\widehat{x}=(\psi,0,\widehat{\alpha_I},\widehat{\alpha_J})^T\in X_0$ be given. Then consider the problem $$\left(\lambda-\left(A_{X_0}+DF_0\left(x_E^0\right)\right))\right)(\varphi,0,\alpha_I,\alpha_J)^T=\widehat x.$$
It re-writes as
\begin{equation*}
\begin{cases}
\lambda \varphi+\varphi^{\prime}+(\mu+\beta_{I}I_{E})\varphi+\beta_{I}s_{E}(.)\alpha_I=\psi,\\
\varphi(0)=0,\\
\lambda \alpha_I=\beta_{I}I_{E}\int_{0}^{\infty}p_I(a)\varphi(a)da+\widehat{\alpha_I},\\
\lambda \alpha_J=\beta_{I}I_{E}\int_{0}^{\infty}p_J(a)\varphi(a)da+\widehat{\alpha_J},
\end{cases}
\end{equation*}  
$\Leftrightarrow$
\begin{equation*}
\begin{cases}
\varphi(a)=\int_0^a e^{-\left(\lambda+\mu+\beta_{I}I_{E}\right)(a-a')}\left[\psi(a')-\beta_{I}s_{E}(a')\alpha_I\right]da',\\
\Delta(\lambda)\alpha_I=\beta_{I}I_{E}\int_{0}^{\infty}p_I(a)\int_0^a e^{-\left(\lambda+\mu+\beta_{I}I_{E}\right)(a-a')}\psi(a')da'da+\widehat{\alpha_I},\\
\lambda \alpha_J=\beta_{I}I_{E}\int_{0}^{\infty}p_J(a)\varphi(a)da+\widehat{\alpha_J}
\end{cases}
\end{equation*}
As a consequence of the above computations, if $\lambda\in\Omega$ satisfies $\Delta(\lambda)\neq 0$ then $\lambda\in\rho\left(A+DF_0\left(x_E^0\right)\right)$.\\
On the other hand, let $\lambda\in\Omega$ be given such that $\Delta(\lambda)=0$. Since since $\lambda\neq 0$ then one can easily check that
\begin{equation*}
\begin{pmatrix}
\varphi_\lambda\\0\\1\\ \alpha_{J,\lambda}\end{pmatrix}\in {\rm Ker}\;\left(\lambda-\left(A+DF_0\left(x_E^0\right)\right)_0\right),
\end{equation*}
wherein we have set 
\begin{equation*}
\varphi_\lambda(a):=-\beta_{I}e^{-\lambda a} s_E(a)\int_0^a e^{\lambda a'}da'\text{ and }\alpha_{J,\lambda}=\frac{\beta_{I}}{\lambda}I_{E}\int_{0}^{\infty}p_J(a)\varphi_\lambda(a)da
\end{equation*}
As a consequence, one deduces that
\begin{equation*}
\sigma\left(A+DF_0\left(x_E^0\right)\right)\cap\Omega= \{\lambda\in\Omega:\;\Delta(\lambda)=0\},
\end{equation*}
and \eqref{eq-suff} follows from Assumption \ref{ASS3} $(iii)$.
This completes the proof of Lemma \ref{claim1}.

\subsection{Proof of Lemma \ref{claim2}}

The proof of the estimate stated in Lemma \ref{claim2} will follow from Gronwall's inequality.
The proof of \eqref{estimation} requires some additional notation. For each $\kappa>0$ we set $M(\kappa)=\frac{\Lambda}{\mu}+\kappa$.
Recalling Definition \eqref{DEF-Lip} and $\tilde B=\overline{B}_{X_0}\left(0,M(\kappa)\right)\cap X_{0+}$, we introduce the quantities:
\begin{equation*}
N_{F_0}(\kappa):=\|F_0\|_{{\rm Lip},\;\tilde B}
\text{ and }N_{G}(\kappa)=\|G\|_{{\rm Lip},\;\tilde B}\;\;\;\;.
\end{equation*}
As a first step let us prove the following lemma:
\begin{lemma}\label{LE-estimate-first}
There exist some constants $M_{**}>0$ and $\rho>0$ such that for each $\kappa>0$, each $\varepsilon\geq 0$ and each $(x,y)\in \left(\overline{B}_{X_0}(0,\kappa)\cap X_{0+}^2\right)$:
\begin{equation}\label{esti-g}
\|U_\varepsilon(t)x-U_\varepsilon(t)y\|\leq M_{**}\|x-y\|e^{\rho M_{**}\left[N_{F_0}(\kappa)+\varepsilon N_{G}(\kappa)\right]t},\;\;\forall t\geq 0.
\end{equation}
\end{lemma}

\begin{proof}
Let $\kappa>0$ and $\varepsilon\geq 0$ be given.
Then for each $(x,y)\in \left(\overline{B}_{X_0}(0,\kappa)\cap X_{0+}\right)^2$ one has:
\begin{equation}\label{variation}
\begin{split}
U_\varepsilon(t)x-U_\varepsilon(t)y&=T_{A_0}(t)(x-y)\\
&+\lim_{\lambda\to\infty}\int_0^t T_{A_0}(t-s)\lambda R\left(\lambda;A\right)\left[F_\varepsilon\left(U_\varepsilon(s)x\right)-F_\varepsilon\left(U_\varepsilon(s)y\right)\right]ds.
\end{split}
\end{equation}

Now recalling estimate $(ii)$ in Theorem \ref{THEO1}, one gets that for each $s\geq 0$:
$$U_\varepsilon(s)x,\;U_\varepsilon(s)y \in \overline{B}_{X_0}\left(0,M(\kappa)\right)\cap X_{0+}.$$
Hence we infer from \eqref{variation} that for each $t\geq 0$:
\begin{equation*}
\begin{split}
&\|U_\varepsilon(t)x-U_\varepsilon(t)y\|\leq \|T_{A_0}(t)\|\|x-y\|\\
&+\limsup_{\lambda\to\infty} \lambda\left\|R(\lambda;A)\right\|_{\mathcal L(X)}\left[N_{F_0}(\kappa)+\varepsilon N_{G}(\kappa)\right]\int_0^t
\|T_{A_0}(t-s)\|\|U_\varepsilon(s)x-U_\varepsilon(s)y\|ds.
\end{split}
\end{equation*}
Recalling that operator $A$ satisfies the Hille-Yosida property, the estimate \ref{esti-g} stated in Lemma \ref{LE-estimate-first} follows from Gronwall's inequality.
\end{proof}

Now note that since $F_0$ is of class $C^2$ one gets that:
\begin{equation}\label{LIP2}
\lim_{\delta\to 0^+}\left\|H_0\right\|_{{\rm Lip}, \overline{B}_{X_0}(x_E^0,\delta)\cap X_{0+}}=0\text{ with }H_0:=F_0(.)-DF_0\left(x_E^0\right).
\end{equation}
Now let us fix $t_0>0$.
Recalling that $U_\varepsilon(t)x_E^\varepsilon\equiv x_E^\varepsilon$, one infers from Lemma \ref{LE-estimate-first} that there exists a constant $M_{**}>0$ such that for each $\delta\in (0,1)$, for each $x\in \overline{B}_{X_0}\left(x_E^\varepsilon,\delta\right)\cap X_{0+}$ and each $\varepsilon\in [0,1]$:
\begin{equation*}
\|U_\varepsilon(t)x-x_E^\varepsilon\|\leq M_{**}\|x-x_E^\varepsilon\|\leq M_{**}\delta,\;\;\forall t\in [0,t_0].
\end{equation*}
As a consequence, since $x_E^\varepsilon\to x_E^0$ as $\varepsilon\to 0$, there exists $\delta_0>0$ small enough such that for each $\delta\in (0,\delta_0)$ and all $\varepsilon\in [0,\delta]$ and all $x\in \overline{B}_{X_0}\left(x_E^0,\delta\right)\cap X_{0+}$:
\begin{equation}\label{esti-m}
\|U_\varepsilon(t)x-x_E^0\|\leq M_{**}\delta+\sup_{\varepsilon\in [0,\delta]}\left\|x_E^\varepsilon-x_E^0\right\|,\;\;\forall t\in [0,t_0].
\end{equation}
Next let $\delta\in (0,\delta_0)$ be given.
Let $x,y\in \overline{B}_{X_0}\left(x_E^0,\delta\right)\cap X_{0+}$ be given.
Then one has for any $z\in\{x,y\}$
\begin{equation*}
U_\varepsilon(t_0)z-L(t_0)z
=\lim_{\lambda\to \infty}\int_0^t T_{A_0}(t-s)\lambda R(\lambda;A)\left[F_\varepsilon\left(U_\varepsilon(s)z\right)-DF_0\left(x_E^0\right)L(s)z\right]ds.
\end{equation*} 
If we set 
\begin{equation*}
\Delta^\varepsilon(t,x,y)=\left[U_\varepsilon(t)x-L(t)x\right]-\left[U_\varepsilon(t)y-L(t)y\right],
\end{equation*}
recalling the definition $H_0$ in \eqref{LIP2} one obtains
\begin{equation*}
\begin{split}
&\Delta^\varepsilon (t,x,y)=\varepsilon \lim_{\lambda\to \infty}\int_0^t T_{A_0}(t-s)\lambda R(\lambda;A)\left[G\left(U_\varepsilon(s)x\right)-G\left(U_\varepsilon(s)y\right)\right]ds\\
&+\lim_{\lambda\to \infty}\int_0^t T_{A_0}(t-s)\lambda R(\lambda;A)\left[H_0\left(U_\varepsilon(s)x\right)-H_0\left(U_\varepsilon(s)y\right)\right]ds\\
&-\lim_{\lambda\to \infty}\int_0^t T_{A_0}(t-s)\lambda R(\lambda;A)DF_0\left(x_E^0\right)\Delta^\varepsilon(s,x,y)ds.
\end{split}
\end{equation*}
If we set $K(\delta)=M_{**}\delta+\sup_{\varepsilon\in [0,\delta]}\left\|x_E^\varepsilon-x_E^0\right\|$, one obtains due to Lemma \ref{LE-estimate-first} and \eqref{esti-m} that there exists some constant $C>0$ (depending on $t_0$ and $\delta_0$) such that for each $t\in [0,t_0]$, each $\varepsilon\in [0,\delta]$:
\begin{equation*}
\|\Delta^\varepsilon (t,x,y)\|\leq C\|x-y\|\left[\varepsilon+m(\delta)\right]+C\int_0^t \|\Delta^\varepsilon(s,x,y)\|ds,
\end{equation*}
wherein we have set
$m(\delta):=\left\|H_0(.)\right\|_{{\rm Lip}, \overline{B}_{X_0}(x_E^0,K(\delta))\cap X_{0+}}$.
Hence Gronwall's inequality applies and provides that for some constant $\tilde C$:
\begin{equation*}
\|\Delta^\varepsilon (t_0,x,y)\|\leq \tilde C \left[\varepsilon+m(\delta)\right]\|x-y\|,
\end{equation*}
for all $x,y\in \overline{B}_{X_0}\left(x_E^0,\delta\right)\cap X_{0+}$, for all $\delta\in (0,\delta_0]$ and for all $\varepsilon\in [0,\delta]$.
As a consequence, one obtains that for all $\delta$ small enough
\begin{equation*}
\sup_{\varepsilon \in [0,\delta]}\left\|U_\varepsilon(t_0)-L(t_0)\right\|_{{\rm Lip},\overline{B}\left(x_E^0,\delta\right)\cap X_{0+}}\leq \tilde C \left[\delta+m(\delta)\right].
\end{equation*}
Recalling that $K(\delta)\to 0$ as $\delta\to 0^+$ and \eqref{LIP2}, we conclude that $m(\delta)\to 0$ as $\delta\to 0^+$ and the result follows.

\section{Application to the prototypical function \eqref{prot1}}  

In this section we will check that Assumption \ref{ASS3} for the prototypical function \eqref{prot1} is indeed satisfied.
We will check $(ii)$ before going to $(iii)$.

\subsection{Assumption \ref{ASS3} $(ii)$}

In this subsection we will derive some conditions on the parameters in order for Assumption \ref{ASS3} $(ii)$ to holds true with the prototypical function \eqref{prot1}.
To do this we first give the following general result:
\begin{proposition}\label{PROP-stability}
In addition to $R_0[0]>1$, we assume that function $p_I\in L_{+}^{\infty}(0,\infty)\cap W^{1,\infty}(0,\infty)$ 
satisfies
\begin{equation}\label{ASS4.1-i}
p_I^{\prime}(a)<\left(  \mu+\beta_{I}I_{E}\right)  p_I(a)\;\;a.e..
\end{equation}
Then Assumption \ref{ASS3} $(ii)$ holds true.
\end{proposition}

Before proving this proposition, let us re-write the hypothesis \eqref{ASS4.1-i} with the prototypical function \eqref{prot1}.
With such a function, $p_I$ takes the form
\begin{equation}\label{prototype}
p_I(a)= 1-\kappa e^{-r a},\;\;\forall a\geq  0. 
\end{equation}
for some $\kappa\in (0,1)$ and $r>0$ .
Together with this function, note that $R_0=R_0[0]$ reads as
\begin{equation*}
R_0=\Lambda\frac{\beta_I}{\mu \nu_I}\left[1-\kappa\frac{\mu }{\mu+r}\right].
\end{equation*}
When $R_0>1$, recall that $\lambda_E$ is uniquely defined by the resolution of the equation 
\begin{equation*}
1=\Lambda\frac{\beta_I}{\nu_I(\mu+\lambda_E)}\left[1-\kappa\frac{\mu+\lambda_E }{\lambda_E+\mu+r}\right].
\end{equation*}
Hence \eqref{ASS4.1-i} re-writes as the following inequality
\begin{equation*}
\left[r+(\mu+\lambda_E)\right]\kappa e^{-r a}< (\mu+\lambda_E),\;\forall a>0.
\end{equation*}
It is then enough to satisfy
\begin{equation}
\frac{r\kappa}{1-\kappa}< (\mu+\lambda_E).
\end{equation}
This condition re-writes as
\begin{equation*}
\Lambda \frac{\beta_I}{\nu_I}p_I^*\left[e^{-\frac{r\kappa}{1-\kappa}.}\right]<1<R_0.
\end{equation*}
The former condition means that using function \eqref{prot1}, Assumption \ref{ASS3} $(ii)$ is satisfied when the endemic value $\lambda_E$ is sufficiently large.

It now remains to prove Proposition \ref{PROP-stability}.\\
\begin{proof}[Proof of Proposition \ref{PROP-stability}]
The proof of this result follows some ideas developed by Magal et al in \cite{MMcCW}.
In this proof, for notational simplicity, we shall denote
\begin{equation*}
x_E^0=\left(s_E,0,I_E,J_E\right).
\end{equation*}
Now let us first notice that coupling bounded dissipavity and asymptotic smoothness in Theorem \ref{THEO1} together with uniform persistence for $U_0$ in Lemma \ref{LE-crucial} $(iii)$, the results in \cite{Magal-Zhao} ensures the existence of an interior global attractor for $U_0$, denoted by $\mathcal A_{00}\subset M_0^0$ (see \eqref{defM00} for the definition of $M_0^0$).

Now we shall construct a suitable Lyapunov functional on $\mathcal A_{00}$ in order to prove that $\mathcal A_{00}=\left\{x_E^0\right\}$.
To that aim let us consider the function $g:x\mapsto x-\ln(x)-1$.
Let $x\in \mathcal A_{00}$ be given and let us consider $\{\left(s(t,.),0,I(t),J(t)\right)^T\}_{t\in\R}\subset \mathcal A_{00}$ an entire solution of the semiflow $\left\{U_0\right\}_{\left\{t\geq 0\right\}}$ such that $\left(s(0,.),0,I(0),J(0)\right)^T=x$.
Next let us consider the uniformly bounded map
\begin{equation}
V[s,0,I,J](t)=\int_{0}^{\infty}\alpha(a)g\left(  \frac{s(t,a)}{s_{E}%
(a)}\right)  da+I_E g\left(  \frac{I(t)}{I_{E}}\right)  ,
\end{equation}
wherein we have set
\begin{equation*}
\alpha(a)\equiv p_I(a)s_{E}(a).
\end{equation*}
Now let us notice that one has
\begin{equation}
\frac{\partial_{t}s}{s_{E}}\left[  1-\frac{s_{E}}{s}\right]  =-\partial
_{a}g\left(  \frac{s}{s_{E}}\right)  +\beta_{I}\frac{s-s_{E}}{s_{E}}\left[
I_{E}-I\right],\;\forall t\in\R.
\end{equation}
Next this yields the following: for each $t\in\R$:
\begin{equation}\label{lyap}
\frac{dV}{dt}\left[s,0,I,J\right](t)= -\int_{0}^{\infty}\alpha(a)\partial_{a}g\left(  \frac{s}{s_{E}}\right)da=\int_{0}^{\infty}\alpha'(a)g\left(  \frac{s}{s_{E}}\right)da.
\end{equation}
Now note that \eqref{ASS4.1-i} re-writes as $\alpha'(a)<0$ for almost every $a>0$.
To complete the proof of Proposition \ref{PROP-stability}, let us consider an increasing sequence $\{t_n\}_{n\geq 0}$ such that $t_n\to-\infty$ as $n\to\infty$. Consider also the sequence of shifted maps $\left(s_n,0,I_n,J_n\right)(t):=\left(s,0,I,J\right)(t+t_n)$. Then possibly along a subsequence, one may assume that
\begin{equation*}
\left(s_n,0,I_n,J_n\right)(t)\to\left(\widehat s,0,\widehat I,\widehat J\right)(t)\text{ as $n\to\infty$}, 
\end{equation*}
locally uniformly with respect to $t\in\R$ with value in $X_{0+}$.
Moreover due to \eqref{lyap} and $\alpha'<0$, one concludes that
\begin{equation*}
\int_{0}^{\infty}\alpha'(a)g\left(  \frac{\widehat s(t,a)}{s_{E}(a)}\right)da=0,\;\;\forall t\in\R.
\end{equation*}
We infer from the former equality that 
\begin{equation*}
\widehat s(t,.)\equiv s_E,\;\;\forall t\in\R.
\end{equation*}
Furthermore function $\left(\widehat I,\widehat J\right)$ satisfies:
\begin{equation*}
\begin{cases}
\left(\partial_a+\mu\right)s_E(a)=-\beta_I\widehat I(t)s_E(a),\;\;t\in\R,\;a>0,\\
\widehat{I}'(t)=\widehat{I}(t)\left(\beta_I p_I^*\left[s_E(.)\right]-\nu_I\right)=0,\\
\widehat{J}'(t)=\beta_I \widehat{I}(t)p_J^*\left[s_E(.)\right]-\nu_J \widehat{J}(t),\;\;t\in\R.
\end{cases}
\end{equation*}
As a consequence, one obtains that
\begin{equation*}
\left(\widehat{I},\widehat{J}\right)(t)=\left(I_E,J_E\right),\;\forall t\in\R.
\end{equation*}
In addition, this yields
\begin{equation*}
\lim_{t\to-\infty}V\left[s,0,I,J\right](t)=V\left[s_E,0,I_E,J_E\right]=0.
\end{equation*}
Finally, since $V\left[s,0,I,J\right](t)\geq 0$ and since this map $t\mapsto V\left[s,0,I,J\right](t)$ is decreasing, one concludes that
\begin{equation*}
V\left[s,0,I,J\right](t)\equiv 0,\;\forall t\in\R,
\end{equation*}
so that $\left(s,0,I,J\right)(t)\equiv \left(s_E,0,I_E,J_E\right)$ for all $t\in\R$ and $x=x_E^0$.
This shows that $\mathcal A_{00}\subset \left\{x_E^0\right\}$ and completes the proof of Proposition \ref{PROP-stability}.
\end{proof}

\subsection{Assumption \ref{ASS3} $(iii)$}
In this subsection we shall discuss Assumption \ref{ASS3}-(iii) with a prototype function \eqref{prot1}, namely $p_I:a\mapsto 1-\kappa e^{-ra}$. In such a case
function $\Delta:\Omega\to \mathbb C$ (see Assumption \ref{ASS3}) reads as
\begin{equation*}
\Delta(\lambda)=\lambda+\beta_{I}^{2}I_{E}\int_0^\infty p_I(a)e^{-(\mu+\beta_{I}I_{E})a}\frac{1-e^{-\lambda a}}{\lambda} da,\;\forall \lambda\in\Omega.
\end{equation*}
Then our next lemma reads as:
\begin{lemma}
For each $\kappa\in [0,1]$ the equation 
\begin{equation*}
\lambda\in\Omega\text{ and }\Delta(\lambda)=0,
\end{equation*}
only has roots with strictly negative real parts.
\end{lemma}

\begin{remark}
One can see $\Delta(\lambda)$ by extension as:
\begin{equation*}
\Delta(\lambda)=\lambda+\beta_{I}^{2}I_{E}\int_0^\infty p_I(a)e^{-(\mu+\beta_{I}I_{E})a}\int_0^a e^{-\lambda a'}da' da,\;\forall \lambda\in\Omega.
\end{equation*}
\end{remark}
\begin{proof}
Let us first notice that
there exists $\varrho>0$ and $M>0$ such that for all $\kappa\in [0,1]$:
\begin{equation*}
\left(\lambda\in\mathbb C\text{ and }\Delta(\lambda)=0\right)\;\Rightarrow\;\left(\varrho\leq |\lambda|\leq M\right).
\end{equation*}
Now note that for $\kappa=0$ the equation only has roots with strictly negative real parts.
Indeed for $\kappa=0$ the equation reads as
\begin{equation*}
\Delta(\lambda)=\lambda+\frac{I_{E}\beta_{I}^{2}\Lambda }{\lambda}\left[\frac{1}{\mu+\lambda_E}-\frac{1}{\mu+\lambda+\lambda_E}\right].
\end{equation*}
Hence recalling that $\lambda_E=I_{E}\beta_{I}$, the equation $\Delta(\lambda)=0$ reduces to
\begin{equation*}
\left(\mu+\lambda_E\right)\lambda^2+\left(\mu+\lambda_E\right)^2\lambda+\beta_{I}\Lambda\lambda_E=0.
\end{equation*}
Thus for $\kappa=0$ the characteristic equation only has roots with strictly negative real parts.

To complete the proof of the lemma we apply Rouch\'e's theorem. It is therefore sufficient to show that for each $\kappa\in [0,1]$ the equation $\Delta(\lambda)=0$ does not have any root on the imaginary axis. \\
Let $\kappa\in [0,1]$ be given. Let us argue by contradiction by assuming that there exists $\omega>0$ such that $\Delta(i\omega)=0$.
Then note that this re-writes as:
\begin{equation*}
-\omega^2+I_{E}\beta_{I}^{2}\Lambda\left[\frac{1}{\mu+\lambda_E}-\frac{1}{\mu+i\omega+\lambda_E}-\frac{\kappa}{r+\mu+\lambda_{E}}+\frac{\kappa}{r+\mu+i\omega+\lambda_{E}} \right]=0.
\end{equation*} 
Taking the imaginary part and since $\omega>0$ yields:
\begin{equation*}
\frac{1}{(\mu+\lambda_E)^2+\omega^2}-\frac{\kappa}{(r+\mu+\lambda_{E})^2+\omega^2}=0.
\end{equation*}
Now the above equality is impossible since $\kappa\in [0,1]$ and $r>0$.
This completes the proof of the lemma.
\end{proof}

\section{Concluding remarks}
In this work, from two systems with age of infection and/or chronological age, we derived an age structured epidemic system with additional mortalities and aggregated infectives modelling the dynamics of transmission of a disease like Hepatitis B virus. The classical forward bifurcation is retrieved: 
\begin{itemize}
\item [(i)] If $R_0\leq 1$, then the system with aggregated infectives has a unique global asymptotically stable disease free equilibrium. Then the disease dies out.\\
\item [(ii)]If $R_0>1$ and we assume that chronic carriers (most of time asymptomatic) have a very low infectious rate, the system with aggregated infectives has (with strong uniform persistence) an unstable disease free equilibrium and a global asymptotically stable endemic equilibrium. 
\end{itemize}

The proof of (ii) uses the useful and complex Theorem 1.2 derived by Magal in \cite{jdemagal}: it shows that it is intuitive but technically difficult to claim that the global stability property of a steady state remains after a slightly perturbation of a parameter.

The result (ii) was an open problem presented in the discrete age formulation of the model studied here\cite[Fig. 3, page 62]{bonzi}. Consider a continuous age model provides a more realistic study of the long term behaviour of the hepatitis B disease\cite{k1}. Moreover the uniform strong persistence of the semiflow is obtained. Roughly speaking here for $R_0 >1$, uniform persistence is the notion saying that the closed subset of extinction for the populations of (acute and chronic) infectives is repelling for the dynamics on the complementary set\cite{jdemagal}. In some sense the persistence points out then the long term survival of infectives over sufficiently large time. 

\section*{Acknowledgements}
The author would like to thank Prof D. B\'{e}koll\`{e} with Prof Houpa D. D. E. for their helpful suggestions which have greatly improved the manuscript. The author address a special acknowledgement to two anonymous faculties for their thorough and independent technical support or critics on an obsolete version of the manuscript. The author is solely responsible for the views and opinions expressed in this research; it does not necessarily reflect the ideas and/or opinions of the funding agencies (AIMS-NEI or IDRC) and University of Ngaoundere.

\end{document}